\documentclass[reqno, 11pt, a4paper]{amsart}
\usepackage{amssymb,graphicx,bbm,color}

\usepackage{ifpdf}
\usepackage[pdftex]{thumbpdf}       

\ifpdf
\usepackage[pdftex,colorlinks]{hyperref}
\hypersetup{%
pdftitle={Exact solution to a Lindley-type equation on a bounded support},
pdfauthor={M. Vlasiou, I.J.B.F. Adan},
pdfkeywords={},
bookmarksnumbered,
pdfstartview={FitH},
urlcolor=cyan,
}%
\else
  \newcommand\phantomsection\relax
\fi

\theoremstyle{plain}
\newtheorem{theorem}{Theorem}

\theoremstyle{remark}
\newtheorem{remark}{Remark}

\newcommand{\p}{\mathbb{P}}

\newcommand{\Dfb}[1][]{\mbox{$F_B^{#1}$}}
\newcommand{\Dfa}[1][]{\mbox{$F_A^{#1}$}}

\newcommand{\Dfw}[1][]{\mbox{$F_W^{#1}$}}
\newcommand{\Dfx}[1][]{\mbox{$F_X^{#1}$}}

\newcommand{\dfw}[1][]{\mbox{$f_W^{#1}$}}

\newcommand{\ltw}{\mbox{$\omega$}}

\newcommand{\eDfb}[1][]{\mbox{$\widehat{F}_B^{#1}$}}
\newcommand{\eDfw}[1][]{\mbox{$\widehat{F}_W^{#1}$}}
\newcommand{\edfw}[1][]{\mbox{$\widehat{f}_W^{#1}$}}

\begin{document}
\title[Lindley-type equation on a bounded support]{Exact solution to a Lindley-type equation on a bounded support}
\author[M.\ Vlasiou, I.J.B.F.\ Adan]{}
\date{September 8, 2005}
\maketitle

{\sc M.\ Vlasiou}\footnote{Corresponding author. Address: LG--1.06, P.O.\ Box 513, 5600 MB Eindhoven, The Netherlands. Email: vlasiou@eurandom.tue.nl}, EURANDOM

{\sc I.J.B.F.\ Adan, Eindhoven University of Technology}

\begin{abstract}
We derive the limiting waiting-time distribution $\Dfw$ of a model described by the Lindley-type equation $W=\max\{0, B - A - W\}$, where $B$ has a polynomial distribution. This exact solution is applied to derive approximations of $\Dfw$ when $B$ is generally distributed on a finite support. We provide error bounds for these approximations.
\end{abstract}

\section{Introduction}\label{s:intro}
Consider a server alternating between two service points. At each service point there is an infinite queue of customers waiting to be served. Only one customer can occupy each service point. Once a customer enters the service point, his total service is divided into two separate phases. First there is a preparation phase, where the server is not involved at all. After the preparation phase is completed the customer is allowed to start with the second phase, which is the actual service. The customer either has to wait for the server to return from the other service point, where he may be still busy with the previous customer, or he may commence with his actual service immediately after completing his preparation phase. This would be the case only if the server had completed serving the previous customer and was waiting for this customer to complete his preparation phase. The server is obliged to alternate; therefore he serves all odd-numbered customers at one service point and all even-numbered customers at the other. Once the service is completed, a new customer immediately enters the empty service point and starts his preparation phase without any delay. In the above setting, the steady-state waiting time of the server $W$ is given by the Lindley-type equation (see also \cite{vlasiou05})
\begin{equation}\label{eq:TheEquation}
W=\max\{0, B - A - W\},
\end{equation}
where $B$ and $A$ are the steady-state preparation and service time respectively.

It is interesting to note  that this equation is very similar to Lindley's equation. The only difference between the two equations is the sign of $W$ at the right hand side. Lindley's equation describes the relation between the waiting time of a customer $W$ and the interarrival time $A$ and service time $B$ in a single server queue. It is one of the fundamental and most well-studied equations in queuing theory. For a detailed study of Lindley's equation we refer to \cite{asmussen-APQ,cohen-SSQ} and the references therein.

The model described by \eqref{eq:TheEquation} applies in many real-life situations that involve a single server alternating between two stations. It was first introduced in \cite{park03}, who study a two-carousel bi-directional system that is operated by a single picker. In this setting, the preparation time $B$ represents the rotation time of the carousels and $A$ is the time needed to pick an item. It is assumed that $B$ is uniformly distributed, while the pick time $A$ is either exponential or deterministic. The authors are mainly interested in the steady-state waiting time of the picker. This problem is further investigated in \cite{vlasiou04}, where the authors expand the results in \cite{park03} by allowing the pick times to follow a phase-type distribution.

In general, it is not possible to derive a closed-form expression for the distribution of $W$ for every given distribution $\Dfb$ of $B$ (or $\Dfa$ of $A$). In \cite{vlasiou05} the authors derive an exact solution under the assumption that $A$ is generally distributed and $B$ follows a phase-type distribution. For the classic Lindley-equation, the M/G/1 single server queue is perhaps the most easy case to analyse. The analogous scenario for our model would be to allow the service time $A$ to be exponentially distributed and the preparation time $B$ to follow a general distribution. For this model though, the analysis is not straightforward, as is the case for Lindley's equation. The structure of $\Dfb$ (or the lack thereof) is essential for this model. If $\Dfb$ belongs to a specific class of distributions, exact computations are possible. This class of distributions includes at least all distribution functions that have a rational Laplace transform and a density on an unbounded support. Both this class and the closed-form expression for the distribution of $W$ are described in detail in \cite{vlasiou05a}.

Despite the fact that this class is fairly big, it does not include all distribution functions. For example, if $\Dfb$ is a Pareto distribution, the method described in \cite{vlasiou05a} is inapplicable. Polynomial distributions are another example of distributions that do not belong to this class. However, they are extremely useful, since they can be used to approximate any distribution function that has a bounded support. Our main goal in this paper is to complement the above mentioned results by deriving a closed-form expression of the steady-state distribution of the waiting time, $\Dfw$,  under the assumption that $A$ is exponentially distributed and $B$ follows a polynomial distribution. In Section \ref{s:Dfw} we derive $\Dfw$ under these assumptions. As an application, in Section \ref{s:Dfw appr} we discuss how one can use this result in order to derive good approximate solutions for $\Dfw$ when $B$ is generally distributed on a bounded support, and we provide error bounds of these approximations. We conclude in Section~\ref{4s:numerical results} with some numerical results.

\section{Exact solution of the waiting time distribution}\label{s:Dfw}
In this section we derive a closed-form expression of $\Dfw$,  under the assumption that $A$ is exponentially distributed and $B$ follows a polynomial distribution. Without loss of generality we can assume that $\Dfb$ has all its mass on $[0,1]$. Therefore, let
\begin{equation}\label{eq:polynomial distr}
\Dfa(x)=1-e^{-\mu x}\qquad\mbox{and}\qquad
\Dfb(x)=\begin{cases}
            \sum_{i=0}^n c_i x^i, &\text{for $0\leqslant x\leqslant 1$;}\\
            1, &\text{for $x \geqslant 1$},
          \end{cases}
\end{equation}
where $\sum_{i=0}^n c_i=1$. Let $X=B-A$. As we have shown in \cite[Section 4]{vlasiou05a}, the mapping
\begin{equation}\label{mapping}
(\mathcal{T}F)(x)=1-\int_x^\infty F(y-x) d\Dfx(y)
\end{equation}
is a contraction mapping --with the contraction constant equal to $\p[B>A]$-- in the space $\mathcal{L}^\infty([0,\infty))$, i.e., the space of measurable and bounded functions on the real line with the norm
$$
\|F\|= \sup_{x \geqslant 0} |F(x)|.
$$
Furthermore, we have shown that $\Dfw$, provided that $\Dfa$ or $\Dfb$ is continuous, is the unique solution to the fixed-point equation $F=\mathcal{T}F$. Then from \eqref{mapping}, for $F=\Dfw$, we have that
\begin{align*}
\Dfw(x) &=1-\int_x^\infty \Dfw(y-x)d\Dfx(y)\\
        &=1-\p[X-W \geq x]=\p[B-W-A \leq x]\\
        &=\int_0^\infty \int_0^\infty \p[B \leq x+z+y] d\Dfa(z) d\Dfw(y)\\
        &=\pi_0 \int_0^\infty \Dfb(x+z) \mu e^{-\mu z} dz  + \int_{0^+}^\infty \int_0^\infty \Dfb(x+y+z) \mu e^{-\mu z} dz d\Dfw(y),
\end{align*}
where $\pi_0$ is the mass of the distribution at the origin, i.e., $\pi_0=\p[W=0]$. Now, by differentiating with respect to $x$, we have after some rewriting (cf.\ \cite[Section 6]{vlasiou05a}) that
\begin{equation}\label{eq:integral}
\dfw(x)=\mu \Dfw(x) - \mu \pi_0 \Dfb(x) - \mu \int_0^\infty \Dfb(x+y) \dfw(y) dy.
\end{equation}
Since $B$ is defined on $[0,1]$, then from Equation \eqref{eq:TheEquation} it emerges that $W$ is also defined on the same interval. Therefore, the integrand at the right-hand side of \eqref{eq:integral} is nonzero only on $[0,1]$. So, substituting \eqref{eq:polynomial distr} in \eqref{eq:integral}, we obtain for $0\leqslant x\leqslant 1$,
\begin{align}\label{eq:f_unkn_poly}
\dfw(x)&=\mu \Dfw(x)-\mu \pi_0 \sum_{i=0}^n c_i x^i-\mu\int_0^{1-x}\sum_{i=0}^n c_i (x+y)^i \dfw(y)dy-\mu\int_{1-x}^1\dfw(y)dy\\
\notag &=\mu \Dfw(x)-\mu \pi_0 \sum_{i=0}^n c_i x^i-\mu\sum_{i=0}^n\sum_{k=0}^i c_i \binom{i}{k}x^{i-k}\int_0^{1-x}y^k\dfw(y)dy - \mu\int_{1-x}^1\dfw(y)dy.
\end{align}
We know from \cite[Section 3]{vlasiou05a} that \eqref{eq:f_unkn_poly} has a unique solution $\dfw$ and $\pi_0$, provided that they satisfy the normalisation equation
\begin{equation}\label{eq:normalisation eq}
\pi_0=1-\int_0^1 \dfw(x)dx.
\end{equation}
To determine $\dfw$ and $\pi_0$, we shall transform the integral equation \eqref{eq:f_unkn_poly} into a (high order) differential equation for $\dfw$. Let $f^{(i)}$ denote the $i$-th derivative of a function $f$. Then differentiating \eqref{eq:f_unkn_poly} with respect to $x$ yields
\begin{align*}
\dfw[(1)](x)&=\mu \dfw(x)-\mu \pi_0 \sum_{i=1}^n i c_i x^{i-1}-\mu \sum_{i=0}^{n-1}\sum_{k=0}^i c_{i+1}(i+1) \binom{i}{k}x^{i-k}\int_0^{1-x}y^k\dfw(y)dy\\
              &\quad+\mu \sum_{i=0}^n\sum_{k=0}^i c_i \binom{i}{k}x^{i-k}(1-x)^k \dfw(1-x)-\mu \dfw(1-x)\\
              &=\mu \dfw(x)-\mu \pi_0 \sum_{i=1}^n i c_i x^{i-1}-\mu \sum_{i=0}^{n-1}\sum_{k=0}^i c_{i+1}(i+1) \binom{i}{k}x^{i-k}\int_0^{1-x}y^k\dfw(y)dy
\end{align*}
and in general, for $\ell=1,2,\ldots,n$,
\begin{equation}\label{eq:derivatives}
\dfw[(\ell)](x)=a_{\ell}(x)+\sum_{j=0}^{\ell-1}\nu_{n-j}(-1)^{\ell-1-j}\dfw[(\ell-1-j)](1-x),
\end{equation}
where
\begin{align*}
&\ \nu_{n-j}=\mu \sum_{i=0}^{n-j}\frac{(i+j)!}{i!}c_{i+j}\\
&\begin{aligned}
 a_{\ell}(x)&=\mu \dfw[(\ell-1)](x)-\mu \pi_0 \sum_{i=0}^{n-\ell}\frac{(i+\ell)!}{i!}c_{i+\ell}\, x^{i}-\mu (-1)^{\ell-1} \dfw[(\ell-1)](1-x)\\
&\qquad-\mu \sum_{i=0}^{n-\ell}\sum_{k=0}^i c_{i+\ell}\frac{(i+\ell)!}{i!} \binom{i}{k}x^{i-k}\int_0^{1-x}y^k\dfw(y)dy.
\end{aligned}
\end{align*}
From \eqref{eq:derivatives} we have that the $n$-th derivative of $\dfw$ is given by
\begin{align*}
\dfw[(n)](x)&=\mu \dfw[(n-1)](x)-\mu\pi_0 n! c_n-\mu (-1)^{n-1}\dfw[(n-1)](1-x)\\
&\quad-\mu n! c_n\int_0^{1-x}\dfw(y)dy+\sum_{j=0}^{n-1}\nu_{n-j}(-1)^{n-1-j}\dfw[(n-1-j)](1-x)\\
&=\mu \dfw[(n-1)](x)-\mu\pi_0 n! c_n-\mu n! c_n\int_0^{1-x}\dfw(y)dy+\sum_{j=1}^{n-1}\nu_{n-j}(-1)^{n-1-j}\dfw[(n-1-j)](1-x),
\end{align*}
which implies that for $0\leqslant x\leqslant 1$,
\begin{align}\label{eq:de}
\dfw[(n+1)](x)&=\mu \dfw[(n)](x)+\mu n! c_n\dfw(1-x)+\sum_{j=1}^{n-1}\nu_{j}(-1)^{j}\dfw[(j)](1-x) \notag\\
              &=\mu \dfw[(n)](x)+\sum_{j=0}^{n-1}\nu_{j}(-1)^{j}\dfw[(j)](1-x).
\end{align}

Up to this point, we have differentiated Equation \eqref{eq:f_unkn_poly} a total of $n+1$ times. Therefore, we need a total of $n+1$ additional conditions in order to guarantee that any solution to \eqref{eq:de} is also a solution to \eqref{eq:f_unkn_poly}. Since for every value of $x$ in $[0,1]$, Equations \eqref{eq:f_unkn_poly} and \eqref{eq:derivatives} are satisfied, then we can evaluate all these equations for a specific $x$, say $x=0$, which provides us with the $n+1$ initial conditions, for $\ell=1,2,\ldots,n$,
\begin{equation}\label{eq:conditions}
\begin{aligned}
\dfw(0)&=\mu\pi_0-\mu\pi_0\,c_0-\mu\sum_{i=0}^n c_i\int_0^{1}y^i\dfw(y)dy\\
\mbox{and}\quad\dfw[(\ell)](0)&=a_{\ell}(0)+\sum_{j=0}^{\ell-1}\nu_{n-j}(-1)^{\ell-1-j}\dfw[(\ell-1-j)](1).
\end{aligned}
\end{equation}
So we now have that Equation \eqref{eq:de} has a unique solution that satisfies these conditions, along with the normalisation equation \eqref{eq:normalisation eq}.

Equation \eqref{eq:de} is a homogeneous linear differential equation, not of a standard form because of the argument $1-x$ that appears at the right-hand side. Therefore, we need to proceed with caution. Note that the unknown probability $\pi_0$ is not involved in \eqref{eq:de}. We shall solve this equation by transforming it into a differential equation we can handle. To this end, substitute $x$ for $1-x$ in \eqref{eq:de}, to obtain the equation
\begin{equation}\label{eq:de2}
\dfw[(n+1)](1-x)=\mu \dfw[(n)](1-x)+\sum_{j=0}^{n-1}\nu_{j}(-1)^{j}\dfw[(j)](x).
\end{equation}
Equations \eqref{eq:de} and \eqref{eq:de2} form a system of equations. Now let
$$
\mbox{\boldmath $\dfw$}(x)=\left[\begin{array}{c}\dfw(x)\\\dfw(1-x)\end{array}\right], \mbox{ }
\mathbf{A_n}=\left[\begin{array}{cc}1 & 0\\0 & (-1)^n\end{array}\right], \mbox{ and }
\mathbf{J}=\left[\begin{array}{cc}0 & 1\\1 & 0\end{array}\right].
$$
Then the system of equations \eqref{eq:de} and \eqref{eq:de2} can be rewritten as
\begin{equation}\label{eq:de in vector}
\mbox{\boldmath $\dfw[(n+1)]$}(x)=\mu \mathbf{A_{n+1}}\mathbf{A_n}\mbox{\boldmath $\dfw[(n)]$}(x)+\mathbf{A_{n+1}}\mathbf{J}\sum_{i=0}^{n-1}\nu_{i}(-1)^{i}\mathbf{A_i}\mbox{\boldmath $\dfw[(i)]$}(x).
\end{equation}
In order to derive the characteristic equation of \eqref{eq:de in vector}, we work as follows. We look for solutions of the form $\mbox{\boldmath $\xi$} e^{r x}$, where $\mbox{\boldmath $\xi$} = \left[ \begin{array}{c}\zeta\\\theta\end{array} \right]$. Substituting this solution into \eqref{eq:de in vector} and dividing by $ e^{r x}$, we derive the following linear system that determines $\mbox{\boldmath $\xi$}$ and $r$, which is
\begin{equation}\label{eq:system for char eq}
 \begin{aligned}
  \zeta r^{n+1}&=\mu \zeta r^{n}+\sum_{i=0}^{n-1}\nu_i \theta r^i\\
    \theta r^{n+1}&=-\mu \theta r^{n}+\sum_{i=0}^{n-1}\nu_i (-1)^{n+1+i}\zeta r^i.
\end{aligned}
\end{equation}
In order for a nontrivial solution to exist, the determinant of the coefficients of $\zeta$ and $\theta$ should be equal to zero. This yields that
\begin{equation}\label{eq:characteristic eq}
r^{2n}(r^2-\mu^2)+(-1)^{n}\left(\sum_{i=0}^{n-1}\nu_i r^i\right)\left(\sum_{j=0}^{n-1}\nu_j (-r)^j\right)=0,
\end{equation}
which is the characteristic equation of \eqref{eq:de in vector}.

Let us \textit{assume} for the moment that the characteristic equation has only simple roots, and label them $r_1,\ldots,r_{2n+2}$. It is interesting to note here that since \eqref{eq:characteristic eq} is a polynomial in $r^2$, then for every root $r$ of this polynomial $-r$ is also a root. Therefore, we shall order the roots so that for every $i$, $r_i=-r_{2n+3-i}$. By substituting each root into the system \eqref{eq:system for char eq}, we obtain the corresponding vectors $\mbox{\boldmath $\xi_i$}$, $i=1,\ldots,2n+2$. Then \eqref{eq:de in vector} has the $2n+2$ linearly independent solutions $\mbox{\boldmath $\xi_i$} e^{r_i x}$. Thus, the general solution of \eqref{eq:de in vector} is given by
\begin{equation}\label{eq:general sol of vector}
\mbox{\boldmath $\dfw$}(x)=\sum_{i=1}^{2n+2} d_i \mbox{\boldmath $\xi_i$} e^{r_i x},
\end{equation}
where $d_i$ are arbitrary constants.

From \eqref{eq:general sol of vector} we can immediately conclude that the solution to Equation \eqref{eq:de} that we are interested in, is of the form
\begin{equation}\label{eq:general sol of de}
\dfw(x)=\sum_{i=1}^{2n+2} d_i \zeta_i e^{r_i x}.
\end{equation}
However, this is not the general solution to \eqref{eq:de}. It does not follow from the derivation of \eqref{eq:general sol of vector} that, for any choice of the coefficients $d_i$, the linear combination \eqref{eq:general sol of de} will satisfy \eqref{eq:de}, since $\zeta_i e^{r_i x}$ is not a solution to \eqref{eq:de}. Therefore, we substitute \eqref{eq:general sol of de} into \eqref{eq:de}, and by keeping in mind that $r_i=-r_{2n+3-i}$, we have that for every $i=1,\ldots,2n+2$,
\begin{equation}\label{eq:relation of d's}
d_i \zeta_i r_i^n (r_i-\mu)=e^{-r_i} d_{2n+3-i}\zeta_{2n+3-i}\sum_{j=0}^{n-1}\nu_j r_i^j.
\end{equation}
These are in fact only $n+1$ relations between the unknown coefficients, since it can easily be shown by using the characteristic equation \eqref{eq:characteristic eq} that the equations for every $i$ and $2n+3-i$ are identical. Using the relations between the coefficients $d_i$, one can rewrite \eqref{eq:general sol of de} as sum of $n+1$ linearly independent solutions to \eqref{eq:de} as follows
\begin{equation}\label{eq:sol of de WITH particulars}
\dfw(x)=\sum_{i=1}^{n+1} d_i \left(\zeta_i e^{r_i x}+ q_i\ \zeta_{2n+3-i} e^{-r_i x}\right),
\end{equation}
where $q_i$ follows from \eqref{eq:relation of d's} if we solve for $d_{2n+3-i}$. Thus, the general solution to \eqref{eq:de} is given by \eqref{eq:sol of de WITH particulars}. The coefficients $d_i$, for $i=1,\ldots,n+1$, and the probability $\pi_0$ that we still need to determine, follow now from the initial conditions \eqref{eq:conditions} and the normalisation equation \eqref{eq:normalisation eq}. Namely, by substituting \eqref{eq:sol of de WITH particulars} to \eqref{eq:conditions} and \eqref{eq:normalisation eq} we obtain a linear system of $n+2$ equations.

Note that it is not possible to use the same argument in order to determine the coefficients $d_i$ for \textit{any} differential equation of the form \eqref{eq:de}, because of its nonstandard form. Here we heavily rely on the fact that we know beforehand that a unique solution exists. We summarise the above in the following theorem.
\begin{theorem} \label{th:density}
Let $\Dfb$ be a polynomial distribution of the form \eqref{eq:polynomial distr}. Then the waiting time distribution $\Dfw$ has a mass $\pi_0$ at the origin, which is given by
\begin{equation*}
\pi_0 = \p[W=0] = 1-\sum_{i=1}^{2n+2}\frac{d_i \zeta_i}{r_i}(e^{r_i}-1),
\end{equation*}
and has a density $\dfw$ on $[0,1]$, given by
$$
\dfw(x)=\sum_{i=1}^{2n+2}d_i \zeta_i e^{r_ix}.
$$
\end{theorem}

Although the roots $r_i$ and coefficients $d_i$ may be complex-valued, the density and the probability $\pi_0$ that appear in Theorem \ref{th:density} will be nonnegative. This follows from the fact that for every distribution $\Dfb$ of the preparation time, \eqref{eq:integral} has a unique solution which is a distribution. It is also clear that, since the differential equation \eqref{eq:de} has real coefficients, then  each root $r_i$ and coefficient $d_i$ have a companion conjugate root and conjugate coefficient, which implies that the imaginary parts cancel.

\begin{remark}
When \eqref{eq:characteristic eq} has roots with multiplicity greater than one, the analysis proceeds essentially in the same way. For example assume that $r_1=r_2$. Then we first look for two solutions to \eqref{eq:de in vector} of the form $\mbox{\boldmath $\xi$} e^{r_1 x}$. If we find only one (that always exists), then we look for a second solution of the form $(x\mbox{\boldmath $\xi$}+\mbox{\boldmath $\eta$}) e^{r_1 x}$, where $\mbox{\boldmath $\eta$}$ is again a vector. Substituting this solution into \eqref{eq:de in vector}, we obtain a linear system that determines $\mbox{\boldmath $\xi$}$ and $\mbox{\boldmath $\eta$}$. Thus we can obtain the general solution to the differential equation \eqref{eq:de in vector}. From this point on, by following the same method, we can formulate a linear system that determines the coefficients $d_i$ and $\pi_0$, and obtain the solution to \eqref{eq:de}.
\end{remark}

\begin{remark}
Another method to derive the solution to the integral equation \eqref{eq:f_unkn_poly} is through Laplace transforms over a bounded interval. We have illustrated this method in \cite{vlasiou04}. The steps of this method are as follows. By taking the Laplace transform of \eqref{eq:f_unkn_poly} over the interval $[0,1]$ we obtain an expression for the Laplace transform $\ltw$ of $\Dfw$ that involves the terms $\ltw(s)$ and $\ltw(-s)$. By substituting $s$ for $-s$ we form a system of two equations from which we can obtain $\ltw$. This step is equivalent to the method we used here, namely forming a system of differential equations for $\dfw(x)$ and $\dfw(1-x)$. It emerges that
$$
\ltw(s)=\frac{P(s)+e^{-s}Q(s)}{R(s)},
$$
where $P$, $Q$, and $R$ are polynomials in $s$. Using the fact that the transform is an analytic function on the whole complex plane, we can deduce that the previous expression is the Laplace transform over a bounded interval of a mixture of $2n+2$ exponentials. This method is fairly straightforward; it is, however, cumbersome and it does not illustrate the special relation between the exponentials with opposite exponents that appear in the density $\dfw$.
\end{remark}

\section{Approximations of the waiting time distribution}\label{s:Dfw appr}
The result we have obtained in the previous section comes in handy in some cases where it is necessary to resort to approximations of the waiting time distribution. We have already proven in \cite{vlasiou05a} that for any distribution of $A$ and $B$ there exists a unique limiting distribution $\Dfw$ for \eqref{eq:TheEquation}, provided that $\p[B<A]>0$, although we may not be able to compute it. Some distributions of the preparation time are not suitable for deriving a closed-form expression of $\Dfw$.  Furthermore, if $\Dfb$ has a bounded support, then we cannot readily apply previously obtained results. In \cite{vlasiou04} only the case where $\Dfb$ is the uniform distribution is covered, while the method described in \cite{vlasiou05a} is not applicable (since distributions on a bounded support are excluded from the class of distributions that are considered there).

Therefore, one may consider approximating $\Dfb$ in order to be able to compute the distribution of the waiting time, which is our main concern. A reasonable approach is to approximate $\Dfb$ by a phase-type distribution. An important reason is that the class of phase-type distributions is dense; any distribution on $[0,\infty)$ can, in principle, be approximated arbitrarily well by a phase-type distribution (see \cite{schassberger-W}). Furthermore, we have shown in \cite{vlasiou05} that if $\Dfb$ is a phase-type distribution, then we can compute explicitly the waiting time distribution $\Dfw$. Nonetheless, if $\Dfb$ has a bounded support, it is more natural and possibly computationally more efficient to fit a polynomial distribution. In the sequel, we shall discuss  how to fit a polynomial distribution to $\Dfb$.

\subsection{Fitting polynomial distributions}
If $\Dfb$ is a continuous distribution on a bounded support, it is reasonable to choose $\eDfb$ to be a polynomial distribution. The famous {\it Weierstrass approximation theorem} asserts the  possibility of uniform approximation of a continuous, real-valued function on a closed and bounded support by some polynomial. The following theorem is a more precise version of Weierstrass' theorem. It is a special case of the theorem by S.\ Bernstein that is stated in \cite[Section VII.2]{feller-IPTIA}.
\begin{theorem}\label{th:bernstein approx}
If $F$ is a continuous distribution on the closed interval $[0,1]$, then as $n\to\infty$
\begin{equation*}
\hat{F}_n(x)=\sum_{k=0}^n F(k/n)\binom{n}{k}x^k(1-x)^{n-k} \to F(x)
\end{equation*}
uniformly for $x\in[0,1]$. Furthermore, $\hat{F}_n$ is also a distribution.
\end{theorem}
\begin{proof}
Bernstein's theorem states that if $F$ is a continuous function, then it can be approximated uniformly in $x$ with the polynomial $\hat{F}_n$. In other words, for any given $\varepsilon>0$, there is an $N$ independent from $x$, such that for all $n>N$, $|\hat{F}_n(x)-F(x)|<\varepsilon$, for all $x$.

It is simple to show that if the function $F$ is a distribution on $[0,1]$, then the approximation $\hat{F}_n(x)$ is also a distribution, since it is continuous, $0\leqslant \hat{F}_n(x)\leqslant 1$, and, by checking its derivative, we shall show that it is non-decreasing in $x$. It suffices to note that
\begin{align*}
\hat{F}_n'(x)&=\sum_{k=1}^n F(k/n)\binom{n}{k} k x^{k-1}(1-x)^{n-k}-\sum_{k=0}^{n-1} F(k/n)\binom{n}{k}x^k(n-k)(1-x)^{n-k-1}\\
&=\sum_{k=0}^{n-1} x^k(1-x)^{n-k-1}\left[ F((k+1)/{n})\binom{n}{k+1}(k+1)-F(k/{n})\binom{n}{k}(n-k) \right]\\
&=\sum_{k=0}^{n-1} x^k(1-x)^{n-k-1}\frac{n!}{k!(n-k-1)!}\left[F((k+1)/{n})-F(k/{n})\right].
\end{align*}
The expression in the square brackets at the right hand side is positive since $F$ is a distribution, which implies that $F((k+1)/{n})\geqslant F(k/{n})$. Therefore, $\hat{F}_n'(x)\geqslant 0$, for $x\in[0,1]$.
\end{proof}
So, given a continuous distribution $\Dfb$ that has all its mass concentrated on $[0,1]$, one can compute a polynomial distribution $\eDfb$ that approximates $\Dfb$ arbitrarily well by using Theorem \ref{th:bernstein approx}. In this sense, the class of polynomial distributions is dense. Then $\eDfw$ can be computed by using Theorem \ref{th:density}.

Naturally, after having obtained an approximation of $\Dfw$, the first question that follows is to determine how good this approximation actually is. Therefore, we shall obtain an upper bound for the error between the approximated distribution for $W$ and the actual one.

\subsection{Bounding the approximation error}
Error bounds for queueing models have been studied widely. The main question is to define an upper bound of the distance between the distribution in question and its approximation, that depend on the distance between the governing distributions. These bounds are obtained both in terms of weighted metrics (see, e.g., \cite{kalashnikov02}) and non-weighted metrics (see, e.g., \cite{borovkov-SPQT, borovkov-ESSP} and references therein). An important assumption which is often made in these studies is that the recursion under discussion should be non-decreasing in its main argument. Clearly, in the model we discuss here this assumption does not hold. Because of Theorem \ref{th:bernstein approx}, we shall limit ourselves to the uniform norm.

Let $\eDfb$ be an approximation of $\Dfb$ and $\eDfw$ the exact solution that we obtain in that case for the distribution of $W$. Let $\widehat{B}$ be a random variable that is distributed according to $\eDfb$, and let $\widehat{X}=\widehat{B}-A$. Define now the mapping (cf.\ \eqref{mapping})
$$
(\widehat{\mathcal{T}}F)(x)=1-\int_x^\infty F(y-x) d F_{\widehat{X}}(y),
$$
which yields that $\eDfw$ is the solution to $F=\widehat{\mathcal{T}}F$ that can be rewritten in the form (cf.\ (2.3))
$$
(\widehat{\mathcal{T}}F)(x)=\int_0^\infty \int_0^\infty \eDfb(x+z+y) \mu e^{-\mu z} dz\,dF(y).
$$
Then we can prove the following theorem.
\begin{theorem}\label{th:error bound}
Let $\|\Dfb-\eDfb\|=\varepsilon$. Then $\|\Dfw-\eDfw\|\leqslant {\varepsilon}/(1-\p[B>A])$.
\end{theorem}
\begin{proof}
We have that
\begin{align*}
\|\Dfw-\eDfw\|&=\|\mathcal{T}\Dfw-\widehat{\mathcal{T}}\eDfw\|=\|\mathcal{T}\Dfw-\mathcal{T}\eDfw+\mathcal{T}\eDfw-\widehat{\mathcal{T}}\eDfw\|\\
              &\leqslant\|\mathcal{T}\Dfw-\mathcal{T}\eDfw\|+\|\mathcal{T}\eDfw-\widehat{\mathcal{T}}\eDfw\|\leqslant \p[B>A]\|\Dfw-\eDfw\|+\|\mathcal{T}\eDfw-\widehat{\mathcal{T}}\eDfw\|,
\end{align*}
since $\mathcal{T}$ is a contraction mapping with contraction constant $\p[B>A]$. Furthermore,
\begin{align*}
\|\mathcal{T}&\eDfw-\widehat{\mathcal{T}}\eDfw\|\\
&=\sup_{x \geqslant 0}\left|\int_0^\infty \int_0^\infty \Dfb(x+z+y) \mu e^{-\mu z} dz d\eDfw(y)- \int_0^\infty \int_0^\infty \eDfb(x+z+y) \mu e^{-\mu z} dz d\eDfw(y)\right|\\
    &\leqslant \sup_{x \geqslant 0} \int_0^\infty\int_0^\infty \mu e^{-\mu z}\left|\Dfb(x+z+y)-\eDfb(x+z+y)\right|dzd\eDfw(y)\\
    &\leqslant \sup_{x \geqslant 0} \int_0^\infty\int_0^\infty \mu e^{-\mu z}\sup_{x+y+z \geqslant 0}\left|\Dfb(x+z+y)-\eDfb(x+z+y)\right|dzd\eDfw(y)\\
    &=\varepsilon\int_0^\infty\int_0^\infty\mu e^{-\mu z}dzd\eDfw(y)=\varepsilon.
\end{align*}
So $\|\Dfw-\eDfw\|\leqslant \p[B>A]\|\Dfw-\eDfw\|+\varepsilon$, which is what we wanted to prove.
\end{proof}

An important feature of Equation \eqref{eq:TheEquation} that made the calculation of an error bound straightforward is that the distribution of the waiting time is the fixed point of a contraction mapping. Note that this is not a property of Lindley's recursion.

\section{Numerical results}\label{4s:numerical results}
This section is devoted to some numerical results. For a given distribution $\Dfb$ we calculate from Theorem~\ref{th:bernstein approx} three polynomial distributions (of first, fifth, and tenth order) that approximate $\Dfb$, and we plot the resulting densities of the waiting time. The distribution $\Dfb$ considered is the piecewise polynomial distribution
$$
\Dfb(x)=(2 x^2)\mathbbm{1}_{[0 \leqslant x \leqslant 1/2]} + (-2 x^2 + 4 x - 1)\mathbbm{1}_{[1/2 \leqslant x\leqslant 1]}+\mathbbm{1}_{[x\geqslant 1]},
$$
where $\mathbbm{1}_{[S]}$ is the indicator function of the set $S$. This distribution is simply the well-known symmetric {\em triangular} distribution on $[0,1]$. Furthermore, we take $\mu=1$.
\begin{figure}[htbp]
\includegraphics[width=0.49\textwidth]{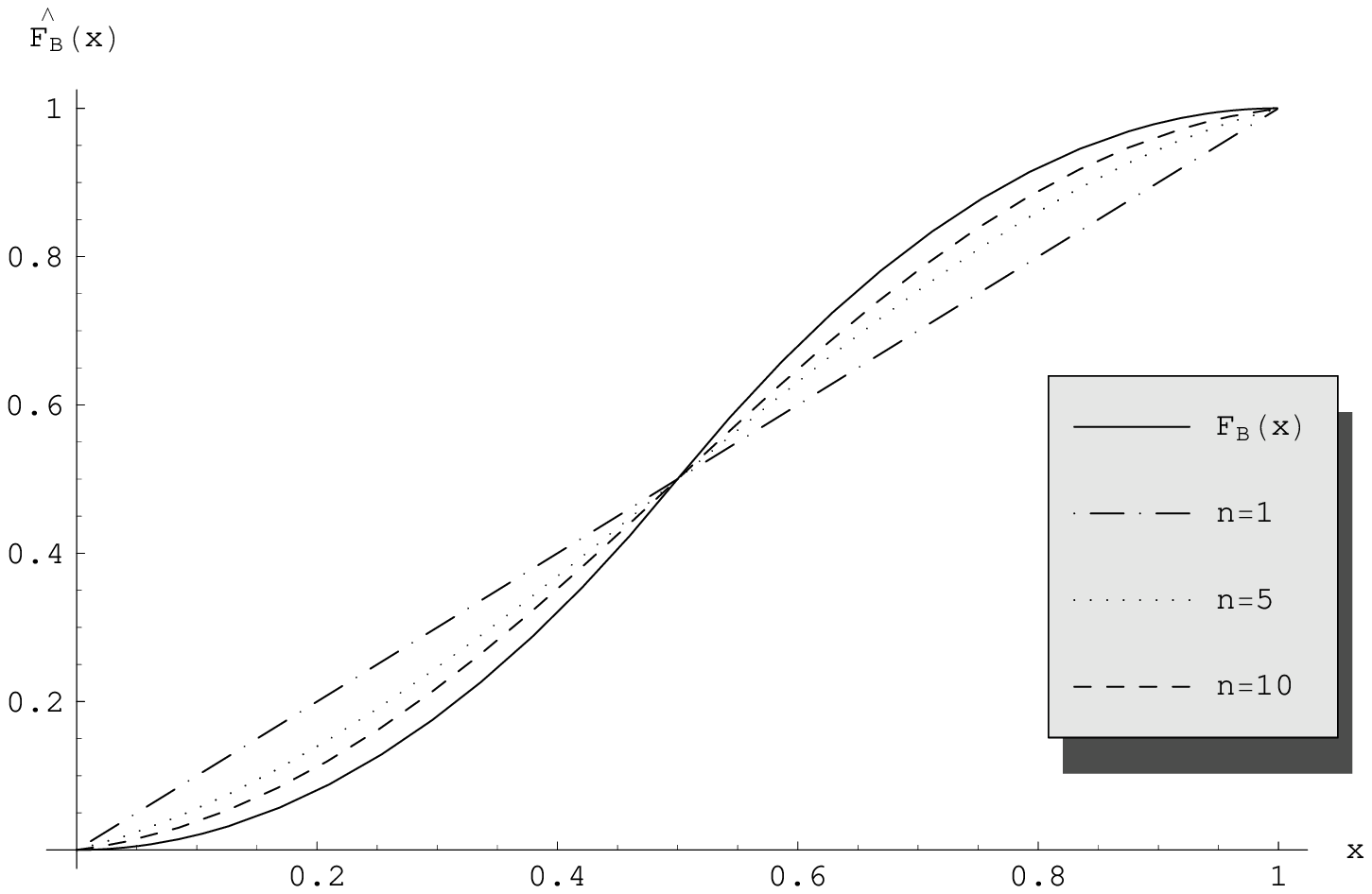}
\hfill
\includegraphics[width=0.49\textwidth]{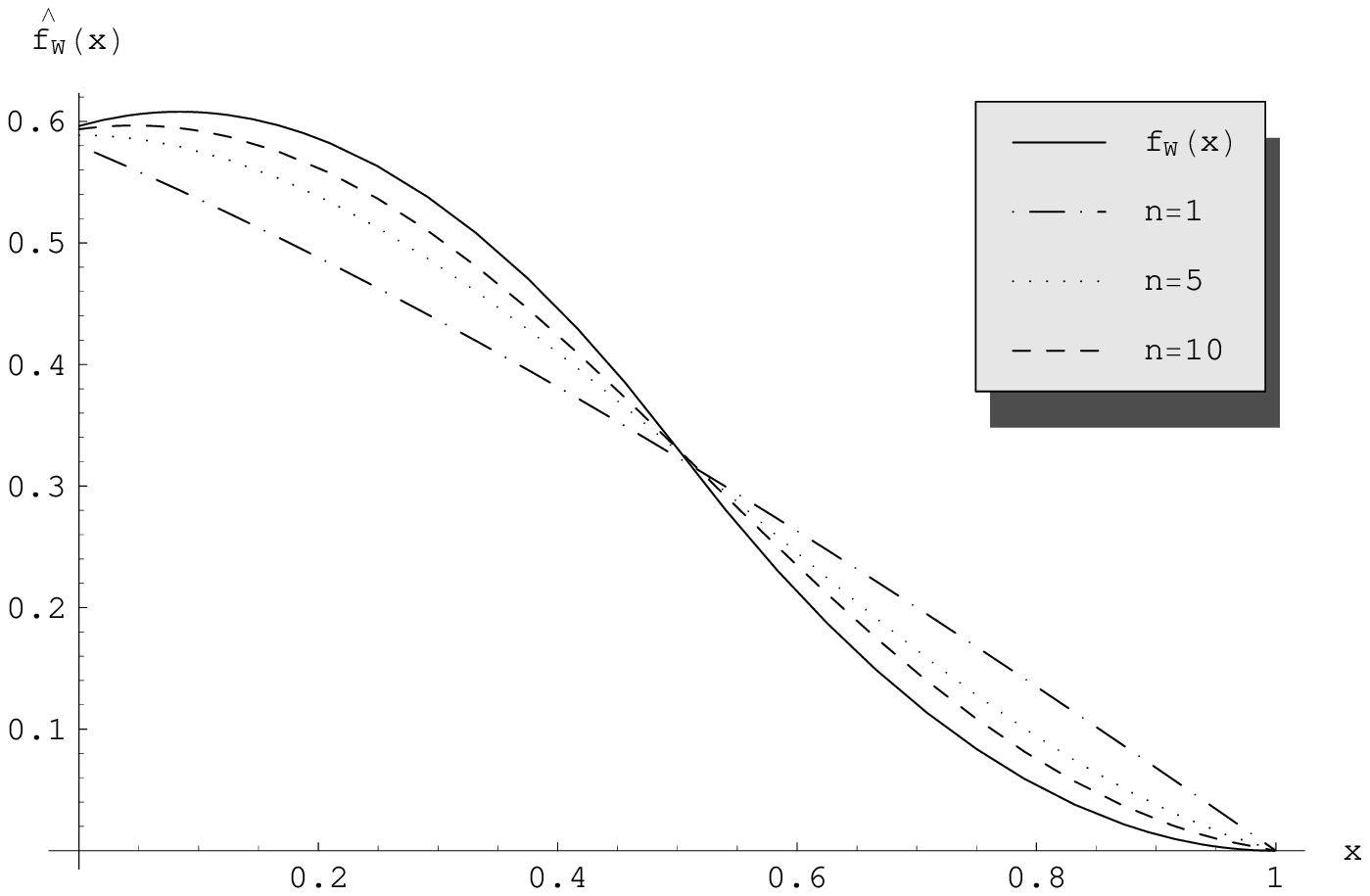}
\caption{The distribution $\Dfb$ and its approximations, and the resulting waiting time densities.}
\label{fig:FB-fW}
\end{figure}

For the above approximations we have computed the distance between $\Dfb$ and $\eDfb$, $\dfw$ and $\edfw$, $\Dfw$ and $\eDfw$, as well as the error bound for $\Dfw$ and $\eDfw$ as it is predicted by Theorem~\ref{th:error bound}. Evidently, the error bound that is predicted by Theorem~\ref{th:error bound} is rather crude. Furthermore, the resulting error between $\Dfw$ and $\eDfw$ in this case is approximately 3 times smaller than the error incurred between the densities (which is an expected consequence of $\pi_0$), and 4.5 times smaller than the initial approximation error between $\Dfb$ and $\eDfb$. Evidently, $A$ smoothes out the error incurred when approximating $\Dfb$. The above are summarised in Table~\ref{table:approximation errors}.
\begin{table}[hbp]
\begin{center}
\renewcommand{\arraystretch}{1.25}
\begin{tabular}{cc|c|c|c}
&\multicolumn{1}{c}{$\|\Dfb-\eDfb\|$}&\multicolumn{1}{c}{$\|\dfw-\edfw\|$}&\multicolumn{1}{c}{$\|\Dfw-\eDfw\|$}&\multicolumn{1}{c}{max $\|\Dfw-\eDfw\|$} \\
$n=1$  &       0.1250       &       0.0841       &       0.0274       &         0.3283         \\
$n=5$  &       0.0664       &       0.0449       &       0.0147       &         0.1744         \\
$n=10$ &       0.0385       &       0.0264       &       0.0086       &         0.1013         \\
\end{tabular}
\end{center}
\caption{The distances between the real distribution or densities and their approximations.\label{table:approximation errors}}
\end{table}

\end{document}